\documentclass{llncs}
\begin{document}

\pagenumbering{arabic}
\pagestyle{headings}
\def\sof{\hfill\rule{2mm}{2mm}}
\def\ls{\leq}
\def\gs{\geq}
\def\SS{\mathcal S}
\def\qq{{\bold q}}
\def\txx{{\frac1{2\sqrt{x}}}}

\title{ {\sc Gamma function, Beta function and combinatorial identities}}
   
\author{Toufik Mansour}
\institute{Department of Mathematics, \\ 
		University of Haifa, Haifa, Israel 31905\\
		{\tt tmansur@study.haifa.ac.il} }
\maketitle
%
\section*{\centering{\sc Abstract}}
In this note we present a method for obtaining a wide class of combinatorial 
identities. We give several examples, in particular, based on the Gamma and Beta 
functions. Some of them have already been considered by previously, and other are new. 

\medskip
\noindent{Key words: Combinatorial identities, Gamma function, Beta function, 
Dirichlet's integral.}
\section*{\centering{\sc 1. Introduction}}
The five basic methods for obtaining combinatorial identities are the following 
(see \cite{ab}). The first is sister Celine's method \cite{FSMC}, generalized by 
Zeilberger ~\cite{Z82}, the second is Gosper's algorithm \cite{G}, Zeilbereger's 
algorithm \cite{Z90,Z91} is the third method, the Wilf and Zeilbereger method \cite{WZ90,WZ92}
is the fourth, and the Hyper algorithm \cite{Pe} is the fifth. Besides, Egorychev \cite{E}
presents another method, based on integral representations in the complex domain. The latter method 
allows to compute several combinatorial sums with inverse binomial coefficients.
In the present note, we present another method for obtaining a wide class 
of combinatorial identities, based on integral representations in the real domain.\\

In 1981, Rockett \cite[Th. 1]{R} (see also \cite{Pl}) proved the following. For any 
nonnegative integer $n$ 
	$$\sum_{k=0}^n {n\choose k}^{-1} = \frac{n+1}{2^{n+1}}\sum_{k=1}^{n+1} \frac{2^k}{k}.$$
In 1999, Trif \cite{T} proved the above result using the Beta function. 
This paper can be regarded as a far-reaching generalization of the ideas presented in \cite{T}. 
Our main result, in its simplest form, can be stated as follows.

\begin{theorem} 
\label{mt01}
Let $r$, $n\geq k$ be any nonnegative integer numbers, let $f(n,k)$ be given by 
	$$f(n,k)=\frac{(n+r)!}{n!} \int_{u_1}^{u_2} p^k(t)q^{n-k}(t) dt,$$
where $p(t)$ and $q(t)$ are two functions defined on $[u_1,u_2]$. 
Let $\{a_n\}_{n\geq 0}$ and $\{b_n\}_{n\geq 0}$ be any two sequences, and let $A(x)$, $B(x)$ 
be the corresponding ordinary generating functions. Then
$$\sum_{n\geq 0} \left[ \sum_{k=0}^n f(n,k) a_kb_{n-k} \right]x^n=
	\frac{d^r}{d x^r} \left[ x^r\int_{u_1}^{u_2} A(xp(t)) B(xq(t)) dt\right].$$
\end{theorem}

As an easy consequence of Theorem \ref{mt01} we get a family of identities, including 
the one presented above.

\begin{example} {\rm (see \cite{JS})}
\label{exa1}
Let $a_n=a^n$ and $b_n=b^n$ for all $n\geq 0$ such that $a+b\neq 0$. 
So the corresponding generating functions of these sequences are 
$$A(x)=\sum\limits_{n\geq 0} a_nx^n=\frac{1}{1-ax};\ \ \ B(x)=\sum\limits_{n\geq 0} b_nx^n=\frac{1}{1-bx}.$$ 

Recall that, the Gamma function is defined by 
		$$\Gamma(s+1)=\int_0^\infty t^se^{-t}dt,$$ 
for all real numbers $s$ such that $s>-1$, and 
the Beta function is defined by 
		$$B(s,r)=\int_0^1 t^{s-1}(1-t)^{r-1}dt,$$ 
for all positive real numbers $s$ and $r$. 
Since $B(s,r)=\frac{\Gamma(s)\Gamma(r)}{\Gamma(s+r+1)}$ú we obtain
	$${s\choose r}^{-1}=(s+1)\int_0^1 t^r(1-t)^{s-r}dt,        \eqno(1)$$
for all nonnegative real numbers $s$ and $r$ such that $s\geq r$.

By Theorem \ref{mt01} and equation $(1)$, 
$$\begin{array}{ll}
\sum\limits_{n\geq0} x^n\sum\limits_{k=0}^n a^kb^{n-k} {n\choose k}^{-1}
	&=\frac{d}{dx} \left( x \int_0^1 \frac{1}{(1-axt)(1-bx+bxt)} dt \right)\\
	&=\frac{d}{dx} \left( \frac{-\ln(1-ax)-\ln(1-bx)}{a+b-abx}  \right).
\end{array}$$
On the other hand, 
	$$\ln(1-x)=\sum\limits_{n\geq 1} \frac{-x^n}{n}\ \mbox{and}\ 
	   \frac{1}{a+b-abx}=\frac{1}{a+b}\sum\limits_{n\geq 0} \frac{a^nb^nx^n}{(a+b)^n},$$ 
hence 
$$\begin{array}{ll}
\displaystyle \sum\limits_{n\geq 0} \sum\limits_{k=0}^n a^kb^{n-k} {n\choose k}^{-1}
	&\displaystyle =\frac{d}{dx} \left[ \frac{1}{a+b} \sum\limits_{n\geq 0} \frac{a^nb^nx^{n+1}}{(a+b)^n}\sum\limits_{k=0}^n \frac{(a^{k+1}+b^{k+1})}{k+1}\cdot\frac{(a+b)^k}{a^kb^k}  \right]\\
	&\displaystyle =\frac{1}{a+b} \sum\limits_{n\geq 0} \frac{(n+1)(ab)^{n+1}x^n}{(a+b)^{n+1}}\sum\limits_{k=1}^{n+1} \frac{(a^k+b^k)(a+b)^{k-1}}{k(ab)^{k-1}},
\end{array}$$
which means that for any nonnegative integer $n$
	$$\sum_{k=0}^n a^kb^{n-k} {n\choose k}^{-1}
		=\frac{n+1}{\left( \frac1a+\frac1b\right)^{n+1}} 
		\sum_{k=1}^{n+1} \frac{(a^k+b^k)\left( \frac1a+\frac1b\right)^{k-1}}{k}.$$
As a numerical example, for $a=b=1$, we get
  $$\sum_{k=0}^n {n\choose k}^{-1}=\frac{n+1}{2^{n+1}}\sum_{k=1}^{n+1} \frac{2^k}{k},$$
which is result of Rockett, Pla and Trif. The latter identity we can 
represent as follows:
  $$\sum_{k=0}^n {n\choose k}^{-1}=(n+1)\sum_{k=0}^n \frac{1}{(n+1-k)2^k}.     \eqno(2)$$
\end{example}

\begin{example}
Let us define $a_n=n$, $b_n=1$ for $n\geq 0$. 
By Theorem \ref{mt01} and equation (1) it is easy to see that
	$$\sum_{n\geq 0} \left[ \sum\limits_{k=0}^n k{n\choose k}^{-1} \right]x^n
			=\frac{-2x\ln(1-x)}{(2-x)^3} -\frac{x(3x-4)}{(2-x)^2(1-x)^2}.$$
Hence, for any nonnegative integer $n$ 
 $$\sum_{k=0}^n k{n\choose k}^{-1}
	=\frac{1}{2^n} \left[ (n+1)(2^n-1)+\sum_{k=0}^{n-2} \frac{(n-k)(n-k-1)2^{k-1}}{k+1} \right].$$ 
\end{example}

In the rest of the paper, we prove Theorem \ref{mt01} and generalize it by functions represented by 
integrals over a $d$-dimensional domain. We present several examples; some of them have been 
considered previously, and other are new. 

\section*{\centering{\sc 2. One-dimensional case}}
First of all, let us prove Theorem \ref{mt01}. Let $f(n,k)$ be as in the statement of the theorem. Then 
 $$\sum_{k=0}^n f(n,k)a_nb_{n-k}
  =\frac{(n+r)!}{n!} \int_{u_1}^{u_2} \sum_{k=0}^n a_kp^k(t) b_{n-k}q^{n-k}(t) dt,$$
which means that
 $$\sum_{n\geq0 } x^n \sum_{k=0}^n f(n,k)a_nb_{n-k}
  =\sum_{n\geq 0} \left[ \frac{(n+r)!x^n}{n!} \int_{u_1}^{u_2} \sum_{k=0}^n a_kp^k(t) b_{n-k}q^{n-k}(t) \right] dt.$$
Let $A(x)=\sum\limits_{n\geq 0} a_nx^n$, $B(x)=\sum\limits_{n\geq 0} b_nx^n$; hence 
$$\sum_{n\geq 0}\sum_{k=0}^n f(n,k) a_kb_{n-k} x^n=
	\frac{d^r}{d x^r} \left[ x^r\int_{u_1}^{u_2} A(xp(t)) B(xq(t)) dt\right],$$
which means that Theorem \ref{mt01} holds.\qed

Now, we present other applications of Theorem \ref{mt01}. 

\begin{example}
\label{ex12}
Immediately, by equation $(1)$ and Theorem \ref{mt01}, we get 
for any nonnegative integer numbers $c$ and $d$  
$$\sum_{n\geq 0} x^{cn}\sum_{k=0}^n {{cn}\choose {dk}}^{-1}=\frac{d}{d x} 
	\int_0^1 \frac{x\cdot dt}{(1-(1-t)^cx^c)(1-t^d(1-t)^{c-d}x^c)}.$$ 
For $c=d=2$ we get 
$$\sum_{n\geq 0} x^{2n}\sum_{k=0}^n {{2n}\choose {2k}}^{-1}=
			\frac{d}{dx}\int_0^1 \frac{x\cdot dt}{(1-(1-t)^2x^2)(1-t^2x^2)}.$$ 
Therefore,
$$\sum_{n\geq 0} x^{2n}\sum_{k=0}^n {{2n}\choose {2k}}^{-1}=
	\frac{d}{dx}\left[ \frac{\ln(1-x)}{x(2-x)}-\frac{\ln(1+x)}{x(2+x)} \right],$$
Hence, for any nonnegative integer $n$ 
$$\sum_{k=0}^n {{2n}\choose {2k}}^{-1}=\frac{2n+1}{2^{2n+2}} \sum_{k=0}^{2n+1} \frac{2^k}{k+1}.$$
\end{example}

\begin{theorem}
\label{mt11}
Let $\{a_n\}_{n\geq 0}$ and $\{b_n\}_{n\geq 0}$ be two sequences, $A(x)$ and $B(x)$ 
be the corresponding ordinary generating functions and $\mu$ be the differential operator of 
the first order defined by $\mu(f)=\frac{d}{dx}(x\cdot f)$. 
Then, for any positive integer $m$
$$\begin{array}{ll}
	\sum\limits_{n\geq 0} & \left[ \sum\limits_{k=0}^n {n\choose k}^{-m} a_kb_{n-k} \right]x^n=\\
		& \mu^m \left[ \underbrace{\int_0^1\int_0^1\cdots\int_0^1}_{m\ \mbox{times}} A(xt_1t_2\dots t_m) B((1-t_1)(1-t_2)\cdots (1-t_m)x) dt_1dt_2\cdots dt_m \right].
\end{array}$$
\end{theorem}
\begin{proof}
Using equation $(1)$ we get
	$${n\choose k}^{-m}=(n+1)^m \left[ \int_0^1 t^k(1-t)^{n-k}dt \right]^{m},$$
which means that 
	$${n\choose k}^{-m}=(n+1)^m \underbrace{\int_0^1\cdots\int_0^1}_{m\ \mbox{times}} (t_1t_2\dots t_m)^k ((1-t_1)(1-t_2)\dots (1-t_m))^{n-k}dt_1\cdots dt_m.$$
So similarly to proof of Theorem \ref{mt01}, this theorem holds.
\qed\end{proof} 

Now let us find another representation for ${n\choose k}^{-m}$. 

\begin{proposition}
\label{mt12}
For any nonnegative integers $n,m$ 
	$$\sum_{k=0}^n {n\choose k}^{-m}=(n+1)^m \sum_{k=0}^n \left[ \sum_{i=0}^k \frac{(-1)^i}{k+1+i} {k\choose i} \right]^m.$$
\end{proposition}
\begin{proof}
By Equation (1) we get for all positive integer $m$ 
	$${n\choose k}^{-m}=(n+1)^m\left( \int_0^1 t^k (1-t)^{n-k} dt \right)^m,$$
which means that
	$${n\choose k}^{-m}=(n+1)^m\left[ \int_0^1 \sum_{i=0}^{n-k} (-1)^i{{n-k}\choose i}t^{k+i}  dt\right]^m,$$
hence the proposition holds.
\qed\end{proof}

The above proposition and equation (2) yield the following.

\begin{corollary}
For any nonnegative integer $n$, 
	$$\sum_{k=0}^n {n\choose k}^{-1}
	=(n+1)\sum_{k=0}^{n} \frac{1}{(n+1-k)2^k}
	=(n+1)\sum_{k=0}^n \sum_{j=0}^k \frac{(-1)^j}{k+1+j}{k\choose j}.$$
\end{corollary}

\begin{corollary}
For any nonnegative integer number $n$,
$$\begin{array}{ll}
\sum\limits_{k=0}^n {n\choose k}^{-2}	&=(n+1)^2\sum\limits_{k=0}^n \left[ \sum\limits_{i=0}^k \frac{(-1)^i}{k+1+i}{k\choose i}\right]^2=\\
				&=(n+1)^2\sum\limits_{k=0}^n \frac{2}{n-k+1}\sum\limits_{j=0}^n \frac{(-1)^j}{n+2+i}{k\choose i}.
\end{array}$$
\end{corollary}
\begin{proof}
By Proposition \ref{mt12} the first equality holds. Now let us prove the 
second equality. By Theorem \ref{mt11} we get  
$$\sum_{n\geq 0} x^n \sum_{k=0}^n {n\choose k}^{-2} =
	\mu^2 \left[ \int_0^1\int_0^1 \frac{1}{(1-tux)(1-(1-t)(1-u)x)} du dt\right],$$
therefore 
      $$\sum_{n\geq 0} x^n \sum_{k=0}^n {n\choose k}^{-2} =
	\mu^2 \left[ \int_0^1 \frac{-2\ln(1-tx)}{x(1-t(1-t)x)} dt\right].$$
Hence, since $\ln(1-tx)=\sum_{n\geq 1} \frac{-t^nx^n}{n}$  
and $\frac{1}{1-t(1-t)x}=\sum_{n\geq 0} t^n(1-t)^nx^n$, the second equality holds.
\qed\end{proof}

As we see, all the above examples depend on equation (1) which follows from the definition 
of Beta function. Now let us present another direction for examples.
Let us define ${s\choose r}$ for any two real numbers $r,s>-1$ as 
	$${s\choose r}=\frac{\Gamma(s+1)}{\Gamma(r+1)\Gamma(s-r+1)}.     \eqno(3)$$

\begin{example}
\label{sin1}
It is easy to see by induction on $n\geq 0$ that 
	$$\int_0^{\pi/2} \sin^{mn}x\,dx
		=\frac{\sqrt{\pi}}{2} \frac{\Gamma((mn+1)/2)}{\Gamma((mn+3)/2)},$$
so Theorem \ref{mt01} for $p(t)=q(t)=a\sin^m t$ yields 
	$$\frac{\sqrt{\pi}}{2} \sum_{n\geq 0} \frac{\Gamma((mn+1)/2)a^n}{\Gamma((mn+3)/2)}=\int_0^{\pi/2} \frac1{1-a\sin^mt} dt.$$
Hence, since $\Gamma(1/2)=\sqrt{\pi}$, we get for any positive integer $m$ 
	$$\sum_{n\geq 0} {{\frac{mn-1}{2}}\choose {\frac{mn}{2}-1}}a^n
	=\frac{2}{\sqrt{\pi}} \int_0^{\pi/2} \frac1{1-a\sin^mt} dt.$$
As numerical examples, we get the following identities:
$$\begin{array}{ll}
	\sum\limits_{n\geq 0} {{\frac{n-1}{2}}\choose {\frac{n}{2}-1}} & =2;\\
	\sum\limits_{n\geq 0} {{n-1/2}\choose {n-1}} & =arctanh(1+\sqrt{2})+arctanh(1-\sqrt{2});\\
	\sum\limits_{n\geq 0} {{\frac{n-1}{2}}\choose {\frac{n}{2}-1}}b^n & =\frac{2\left[ arctanh\left( \frac{b-1}{\sqrt{b^2-1}} \right) -arctanh\left( \frac{b}{\sqrt{b^2-1}} \right) \right] }{\sqrt{b^2-1}},\ \ b>1.
\end{array}$$
\end{example}
\section*{\centering{\sc 3. Generalization: {$d$}-dimensional case}}
Let $X$ be a multiset of variables $x_j$, where $j=1,2,\dots,d+1$, and let 
$X'=\{x_{i_1},\dots,x_{i_l}\}$ be the underlying set. Let $g(t)$ and $f_j(t)$, 
$j=1,2,\dots,d$ be any $d+1$ functions such that 
$\phi(x_{i_1},\dots,x_{i_l})=g(x_{d+1})\prod\limits_{j=1}^d f_j(x_j)$ is a function 
defined on a $l$-dimensional domain $D$. 
Let $r$ be a nonnegative integer number, and let $f(k_1,k_2,\dots,k_d)$ be given by 
        $$f(k_1,k_2,\dots,k_d)=\frac{(k_1+\dots+k_d+r)!}{(k_1+\dots+k_d)!} 
		\int_D \phi(x_{i_1},\dots,x_{i_l})\,dx_{i_1}\dots dx_{i_l}.$$
Then for any sequences $\{a_n^{(j)}\}_{n\geq 0}$, $j=1,2,\dots d$, 
$$\begin{array}{ll}
	\sum\limits_{k_1+\dots+k_d=n} & f(k_1,k_2,\dots,k_d)\prod\limits_{j=1}^d a^{(j)}_{k_j} =\\
		& =\frac{(n+r)!}{n!}\int_D \left( g(x_{d+1})\sum\limits_{k_1+\dots+k_d=n} \prod\limits_{j=1}^d a_{k_j}^{(j)} f_j^{k_j}(x_j) \right) dx_{i_1}\dots dx_{i_l}.
\end{array}$$
Therefore 
$$\begin{array}{ll}
	\sum\limits_{n\geq 0}\sum\limits_{k_1+\dots+k_d=n} & f(k_1,k_2,\dots,k_d)x^n\prod\limits_{j=1}^d a^{(j)}_{k_j} =\\
		&=\sum\limits_{n\geq 0} \left[ \frac{(n+r)!}{n!} \int_D g(x_{d+1})\sum\limits_{k_1+\dots+k_d=n} \prod\limits_{j=1}^d a_{k_j}^{(j)}(xf_j(x_j))^{k_j} dx_{i_1}\dots dx_{i_l} \right]
\end{array},$$
which means that
$$\begin{array}{ll}
	\sum\limits_{n\geq 0}\sum\limits_{k_1+\dots+k_d=n} & f(k_1,k_2,\dots,k_d)x^n\prod_{j=1}^d a^{(j)}_{k_j} =\\
		&=\frac{d^r}{d x^r} \left[ x^r \int_D g(x_{d+1})\prod\limits_{j=1}^d A_j(xf_j(x_j)) dx_{i_1}\dots dx_{i_l} \right] 
\end{array},$$
where $A_j(x)$ is the generating function of the sequence $\{a_n^{(j)}\}_{n\geq 0}$.
Hence, we get the following result (Theorem \ref{mt01} is its particular case), 
which gives us a general method for obtaining combinatorial identities.

\begin{theorem} 
\label{method}
Let $X$ be a multiset of variables $x_j$, where $j=1,2,\dots,d+1$, and let 
$X'=\{x_{i_1},\dots,x_{i_l}\}$ be the underlying set. Let $g(t)$ and $f_j(t)$, 
$j=1,2,\dots,d$ be any $d+1$ functions such that 
$\phi(x_{i_1},\dots,x_{i_l})=g(x_{d+1})\prod\limits_{j=1}^d f_j(x_j)$ is a function 
defined on a $l$-dimensional domain $D$. 
Let $r$ be a nonnegative integer number, and let $f(k_1,k_2,\dots,k_d)$ be given by 
        $$f(k_1,k_2,\dots,k_d)=\frac{(k_1+\dots+k_d+r)!}{(k_1+\dots+k_d)!} 
		\int_D \phi(x_{i_1},\dots,x_{i_l})\,dx_{i_1}\dots dx_{i_l}.$$
Then for any sequences $\{a_n^{(j)}\}_{n\geq 0}$, $j=1,2,\dots d$, 
$$\begin{array}{ll}
\sum\limits_{n\geq 0}\sum\limits_{k_1+\dots+k_d=n} & f(k_1,k_2,\dots,k_d)x^n\prod\limits_{j=1}^d a_{k_i}^{(j)} =\\
		&\frac{d^r}{d x^r} \left[ x^r  \int_D g(x_{d+1})\prod\limits_{j=1}^d A_j(xf_j(x_j)) dx_{i_1}\dots dx_{i_l}\right] ,
\end{array}$$
where $A_j(x)$ is the ordinary generating function of the sequence $\{a_n^{(j)}\}_{n\geq 0}$.
\end{theorem}

Another way to generalize Theorem \ref{mt01} is the following.
Let $V$ be the hyperplane defined by $\sum_{i=1}^d \left( \frac{x_i}{a_i} \right)^{p_i}=1$ 
where $x_i\geq 0$ for all $i=1,2,\dots d$. If $p_i\geq 0$ for all $i$, 
then the {\em Dirichlet's integral} is defined by 
$$\int_V \prod_{j=1}^d x_j^{\alpha_j-1} dx_1\cdots dx_d=\frac{a_1^{\alpha_1}\cdots a_d^{\alpha_d}}{p_1\cdots p_d}\frac{\Gamma\left(\frac{\alpha_1}{p_1}\right)\cdots \Gamma\left(\frac{\alpha_d}{p_d}\right)}{\Gamma\left(1+\frac{\alpha_1}{p_1}+\cdots+\frac{\alpha_d}{p_d}\right)}.      \eqno(4)$$
So for $p_j=1$, $a_j=1$, and $\sum_{j=1}^d \alpha_j=n$ we obtain 
	$${n\choose {\alpha_1,\cdots,\alpha_d}}^{-m}=\frac{(n+d-1)!^m}{n!^m} \left( \int_{x_1+\cdots+x_d=1} x_1^{\alpha_1}\cdots x_d^{\alpha_d}dx_1\cdots dx_d \right)^m.        \eqno(5)$$
Hence, Theorem \ref{method}, Theorem \ref{mt01} and equation (4) yield the following.

\begin{theorem}
\label{geth}
Let $\{a_n^{(j)}\}_{n\geq 0}$ be any sequence for all $j=1,2,\dots,d$, and let $\nu$ be the differential operator 
of the $(d-1)$th order defined by $\nu_d(f)=\frac{d^{d-1}}{dx^{d-1}}(x^{d-1}f)$ . Then 
$$\begin{array}{ll}
\sum\limits_{n\geq 0} x^n & \sum\limits_{\alpha_1+\cdots +\alpha_d=n} {n\choose {\alpha_1,\cdots,\alpha_d}}^{-m} \prod_{j=1}^d a_{\alpha_j}^{(j)}=\\
			  &= \nu_d^m \left[\underbrace{\int_V\cdots\int_V}_{m\ \mbox{times}} 
				\prod\limits_{j=1}^d A_j(xx_{j,1}x_{j,2}\cdots x_{j,m})  \prod_{i=1,j=1}^{d,m} dx_{i,j} \right],
\end{array}$$ 
where $V$ is the hyperplane defined by $x_1+x_2+\dots+x_d=1$ ,  $A_j(x)$ is the ordinary generating function of sequence $\{a_n^{(j)}\}_{n\geq 0}$, $j=1,2,\dots,d$.
\end{theorem}


\begin{thebibliography}{99}
\bibitem[JS]{JS}
{\sc C.H. Jinh, and L.J. Sheng}, On sums of the inverses of binomial coefficients, 
{\em Tamkang J. Management Sci.} {\bf 8} (1987), no. 1, 45--48.

\bibitem[E]{E}
{\sc G.P. Egorychev}, Integral representation and the computation of combinatorial sums, 
{\em Translations of mathematical monographs} {\bf 59}, Amer. Math. Soc. 1984. 

\bibitem[F]{FSMC}
{\sc Fasenmyer, Sister mary Celine}, Some generalized hypergeometric polynomials, 
Ph.D. dissertation, University of Michigan, November, 1945.

\bibitem[G]{G}
{\sc R.W. Gosper}, Decision procedure for indefinite hypergeometric summation, 
{\em Proc. Natl. Acad. Sci. USA} {\bf 75} (1978), 40--42.

\bibitem[Pl]{Pl}
{\sc J. Pla}, The sum of the inverses of binomial coefficients revisited, 
{\em The Fibonacci Quarterly} {\bf 35} (1997) 342--345.

\bibitem[Pe]{Pe}
{\sc M. Petkov\u sek}, Finding closed-form solutions of difference equations 
by symbolic methods, Ph.D. thesis, Carnegie-Mellon University, CMU-CS-91-103, 1991.

\bibitem[PWZ]{ab} 
{\sc M. Petkov\v sek, H.S. Wilf, and D. Zeilbereger}, $A=B$. 
A K Peters, Ltd., Wellesley, MA, 1996. 

\bibitem[R]{R}
{\sc A.M. Rockett}, Sums of the inverses of binomial coefficients, 
{\em The Fibonacci Quarterly} {\bf 19} (1981) 433--437.

\bibitem[T]{T}
{\sc T. Trif}, Combinatorial sums and series involving inverses of binomial coefficients, 
{\em The Fibonacci Quarterly} {\bf 38} (2000) 79--84.

\bibitem[WZ90]{WZ90}
{\sc H. Wilf and D. Zeilberger}, Rational functions certify combinatorial identities,
{\em J. Amer. Math. Soc.} {\bf 3} (1990), 147--158.

\bibitem[WZ92]{WZ92}
{\sc H. Wilf and D. Zeilberger}, An algorithmic proof theory for hypergeometric 
(ordinary and "q") multisum/integral identities, 
{\em Inv. Math.} {\bf 108} (1992), 575--663.

\bibitem[Z82]{Z82}
{\sc D. Zeilberger}, Sister Celine's technique and its generalization, 
{\em J. Math. Anal. Apll.} {\bf 85} (1982), 114--145.

\bibitem[Z90]{Z90}
{\sc D. Zeilberger}, A fast algorithm for proving terminating hypergeometric 
identities, {\em Discr. Math.} {\bf 80} (1990), 207--211.

\bibitem[Z91]{Z91}
{\sc D. Zeilberger}, The method of creative telescoping, 
{\em J. Symbolic Computation} {\bf 11} (1991), 195--204.
\end{thebibliography}
\end{document}